\documentclass{amsproc}

\newtheorem{theorem}{Theorem}[section]
\newtheorem{lemma}[theorem]{Lemma}
\newtheorem{proposition}[theorem]{Proposition}

\theoremstyle{definition}

\newtheorem{assumption}[theorem]{Assumption}

\theoremstyle{remark}
\newtheorem{remark}[theorem]{Remark}
\newtheorem{conjecture}[theorem]{Conjecture}

\newcommand{\R}{{\mathbb R}}

\newcommand{\Z}{{\mathbb Z}}
\newcommand{\FF}{{\mathcal F}}
\newcommand{\GG}{{\mathcal G}}
\newcommand{\HH}{{\mathcal G}}
\newcommand{\Int}{{{\rm Int}\,}}

\newcommand{\N}{{\mathbb N}}

\newcommand{\Fix}{{\rm Fix}}

\newcommand{\BB}{{\mathcal B}}

\newcommand{\DD}{{\mathcal D}}
\newcommand{\DC}{{\rm Diff}_c^r}

\newcommand{\Diff}{{\rm Diff}}
\newcommand{\bd}{{\sc Proof}.\ \ }

\title[Actions of groups of diffeomorphisms on one-manifolds]
{Actions of groups of diffeomorphisms on one-manifolds by $C^1$ diffeomorphisms}

\author[S. Matsumoto]{Shigenori Matsumoto}
\address{Shigenori Matsumoto, 
Department of Mathematics, College of
Science and Technology, Nihon University, 1-8-14 Kanda, Surugadai,
Chiyoda-ku, Tokyo, 101-8308 Japan
}
\email{matsumo@math.cst.nihon-u.ac.jp
}
\thanks{The author is partially supported by Grant-in-Aid for
Scientific Research (C) No.\ 25400096.}


\begin{document}

\begin{abstract}
Denote by $\DC(M)_0$ the identity component of the group
of the compactly supported $C^r$ diffeomorphisms of a
connected $C^\infty$ manifold $M$.
We show that if $\dim(M)\geq2$ and $r\neq \dim(M)+1$, then any homomorphism
from $\DC(M)_0$ to ${\Diff}^1(\R)$ or ${\Diff}^1(S^1)$ is trivial.
\end{abstract}

\maketitle

\section{Introduction}

\'E. Ghys \cite{G} asked if the group of diffeomorphisms of a manifold
admits a nontrivial action on a lower dimensional manifold.
A break through towards this problem was obtained by K. Mann \cite{M}
for one dimensional target manifolds.
Let $M$ be a connected $C^\infty$ manifold without
boundary, compact or not. For $r=0,1,2,\cdots,\infty$,
denote by ${\rm Diff}^r_c(M)_0$ the identity component of the group
of the compactly supported $C^r$  diffeomorphisms (homeomorphisms for
$r=0$) of $M$.

\begin{theorem} \label{tm} {\rm (K. Mann)} \ \
Any homomorphism from ${\rm Diff}^r_c(M)_0$ to
${\rm Diff}^2(S^1)$ or to ${\rm Diff}^2(\R)$ is trivial,
provided $\dim(M)\geq 2$ and $r\neq \dim(M)+1$.
\end{theorem}

For a simpler proof of this fact,
see also \cite{M''}.  
A natural question is whether it is possible to lower the
differentiability
of the target group. 
In fact for $r=0$, E. Militon \cite{Mi2} obtained the final result.

\begin{theorem} {\rm (E. Militon)}
Any homomorphism from ${\rm Diff}^0_c(M)_0$ to
${\rm Diff}^0(S^1)$ is trivial
if $\dim(M)\geq 2$.
\end{theorem}

Notice that ${\rm Diff}^0(\R)$ can be considered to be a subgroup
of ${\rm Diff}^0(S^1)$. So we do not mention in the above theorem
the case where the target group is ${\rm Diff}^0(\R)$.

Even for $r\geq1$, we have:

\begin{conjecture}
Any homomorphism from ${\rm Diff}^r_c(M)_0$ to
${\rm Diff}^0(S^1)$ is trivial
if $\dim(M)\geq 2$. 
\end{conjecture}

The purpose of this paper is to mark one step forward towards this conjecture.

\begin{theorem} \label{T}
If $\dim(M)\geq2$ and $r\neq \dim(M)+1$,
any homomorphism
from $\DC(M)_0$ to $\Diff^1(S^1)$ or $\Diff^1(\R)$ is trivial.
\end{theorem}

Frequent use of the simplicity of the group $\DC(M)_0$ is made in
the proof.
The condition $r\neq \dim(M)+1$ is needed for it.
As for Theorem \ref{tm}, the proof is built upon a theorem of Kopell and 
Szekeres about  $C^2$ actions of abelian groups on a compact interval,
while for Theorem \ref{T}, upon a theorem of Bonatti, Monteverde, Navas
and Rivas about 
$C^1$ actions of solvable  Baumslag-Solitar groups on a compact interval.

By virtue of the fragmentation lemma, Theorem \ref{T} reduces to:

\begin{theorem} \label{T'}
For $n\geq2$ and $r\neq n+1$, any homomorphism from
$\Diff^r_c(\R^n)_0$ to $\Diff^1(S^1)$ or $\Diff^1(\R)$ is trivial.
\end{theorem}

In Section 2, we show that the case of target group $\Diff^1(S^1)$ can
be reduced to the case $\Diff^1(\R)$.
In Sections 3 and 4, we establish fixed point results for certain
subgroups of $\Diff^\infty_c(\R^n)_0$.
In Section 5, we prove Theorem \ref{T'} following
an argument of E. Militon \cite{Mi2}.
Finally we give some sporadic results for $\Diff^0(S^1)$ target
in Section 6.

\medskip
{\sc Acknowledgement:} The author is greatful to Kathryn Mann and
Andres Navas for helpful comments and conversations.

\section{Reduction to the case $\Diff^1(\R)$}

In this section, we show that Theorem \ref{T'} for the target
group $\Diff^1(S^1)$ is reduced to the case of $\Diff^1(\R)$.

\begin{proposition}\label{p0}
Let $r\neq n+1$ and $n\geq1$.
Assume that $\Phi:\Diff^r_c(\R^n)_0\to\Diff^0(S^1)$ is a nontrivial homomorphism.
then the global fixed point set is nonempty:
$\Fix(\Phi(\Diff^r_c(\R^n)_0))\neq\emptyset$.
\end{proposition}

This proposition enables us to conclude that the image of $\Phi$
is contained in the group of
the homeomorphisms of $\R$. In particular,
Theorem \ref{T'} for the target group $\Diff^1(S^1)$
is reduced to the case of $\Diff^1(\R)$.

Denote $\GG=\Diff^r_c(\R^n)_0$.
By the simplicity of the group $\GG$, 
the homomorphism $\Phi$ in the proposition is
injective and its image is contained in $\Diff^0_+(S^1)$,
the group of the orientation preserving homeomorphisms.

Let
$B_0$  be the closed unit ball in $\R^n$ centered at the origin.
Define a family $\BB$ of the closed balls in $\R^n$ by
$$\BB=\{g(B_0)\mid g\in\GG\}.$$
Also for $B\in\BB$, let
$$\GG(B)=\{g\in\GG\mid {\rm Supp}(g)\subset\Int(B)\}.$$

To show Proposition \ref{p0}, it is sufficent to show the following.

\begin{proposition}\label{p10}
For any $B\in\BB$, the fixed point set $\Fix(\Phi(\GG(B)))$ is nonempty.
\end{proposition}

In fact, choose an increasing sequence of balls,
$\{B_k\}_{k\in\N}\subset\BB$
such that $\bigcup_k B_k=\R^n$. Then we have $\GG=\cup_k\GG(B_k)$
and $\Fix(\Phi(\GG))=\bigcap_k\Fix(\Phi(\GG(B_k)))$. Therefore
by the compactness of $S^1$, Propositon \ref{p0} follows from
Proposition \ref{p10}.

Now for any $B_1,B_2\in\BB$, the groups $\GG(B_1)$ and $\GG(B_2)$ are
conjugate in $\GG$. Therefore their images $\Phi(\GG(B_1))$ and
$\Phi(\GG(B_2))$ are conjugate in $\Diff^0_+(S^1)$. 
They are simple. Moreover if $B_1$ and $B_2$ are disjoint, any
element of $\Phi(\GG(B_1))$ commutes with any element of
$\Phi(\GG(B_2))$. Therefore Proposition \ref{p10} reduces to
the following.

\begin{proposition} \label{p1}
Let $G_1$ and $G_2$ be simple nonabelian
subgroups of $\Diff^0_+(S^1)$.
Assume that $G_2$ is conjugate to $G_1$ in
$\Diff^0_+(S^1)$ and that any element of $G_1$ commutes with any element of $G_2$.
Then there is a global fixed point of $G_1$:
$\Fix(G_1)\neq\emptyset$.
\end{proposition}

{\sc Proof}. 
Let $X_2\subset S^1$ be a minimal set of $G_2$.
The set $X_2$ is either a finite set, a Cantor set or
the whole of $S^1$. If $X_2$ is a singleton, then
$G_2$ admits a fixed point. Since $G_1$ is conjugate to
$G_2$, we have $\Fix(G_1)\neq\emptyset$, as is required.
So assume for contradiction that $X_2$ is not a singleton.

First if $X_2$ is a finite set which is not a singleton, we
get a nontrivial homomorphism from $G_2$ to a
finite abelian group, contrary to the assumption of the simplicity.
In the remaining case, it is well known, easy to show,
that the minimal set is unique. That is, $X_2$
is contained in any nonempty $G_2$ invariant closed subset.

Let $F_1$ be the subset of $G_1$ formed by the elements $g$
such that $\Fix(g)\neq\emptyset$.
Let us show that there is a nontrivial element in $F_1$.
Assume the contrary. Then $G_1$ acts freely on $S^1$.
Consider the group $\tilde G_1$ formed
by any lift of any element of $G_1$ to the universal covering space
$\R\to S^1$.
Now $\tilde G_1$ acts freely on $\R$. 
A theorem of H\"older asserts that $\tilde G_1$ is abelian.
See \cite{N} for a short proof, or \cite{Th} for an even shorter proof.
The canonical projection
$\pi:\tilde G_1\to G_1$ is a group homomorphism,
and $G_1=\pi(\tilde G_1)$
would be abelian, contrary to the assumption of the proposition.
 
Since $G_1$ and $G_2$ commutes, the fixed point set
$\Fix(g)$ of any element $g\in F_1$ is $G_2$ invariant.
 Therefore we have
\begin{equation}\label{ee}
\mbox{ $X_2\subset\Fix(g)$ for any $g\in F_1$}.
\end{equation}
This shows that $F_1$ is 
in fact a subgroup. By the very definition,
$F_1$ is normal. Since $G_1$ is simple and $F_1$ is nontrivial, $F_1=G_1$.
Finally again by (\ref{ee}), $\Fix(G_1)\neq\emptyset$. 
Then the minimal set of $G_1$ must be a singleton. Since
$G_2$ is conjugate to $G_1$, the minimal set $X_2$ of $G_2$ 
is also a singleton, contrary to the assumption. \qed

\section{Fixed point set of $\Phi(G)$}
Again consider $\HH=\Diff^r_c(\R^n)_0$, where $n\geq1$ and $r\neq n+1$.
We shall show Theorem \ref{T'} for the target group ${\rm Diff}^1(\R)$
by a contradiction.
So let us assume that 
$\Phi:\HH\to\Diff^1(\R)$ is a nontrivial homomorphism.
By the simplicity of $\HH$, $\Phi$ is injective and its image is contained
in $\Diff^1_+(\R)$.
For the purpose of showing Theorem \ref{T'}, it is no loss of generality to assume the following.

\begin{assumption} \label{a}
There is no global fixed point of $\Phi(\HH)$: $\Fix(\Phi(\HH))=\emptyset$.
\end{assumption}

In fact, we only have to pass from $\R$ to a connected
component of $\R\setminus\Fix(\Phi(\HH))$.
This assumption will be made all the way until the end of the proof of
Theorem \ref{T'}.

We consider an embedding of Baumslag-Solitar group ${\rm BS}(1,2)$ into
the group $\HH(B)$. See Section 2 for the definition of $\HH(B)$.
Recall that

$${\rm BS}(1,2)=\langle a,b\mid aba^{-1}=b^2\rangle.$$

This group is a subgroup of $GA$, the group of the orientation preserving
affine transformations of $\R$, where $a$ corresponds to $x\mapsto 2x$,
and $b$ to $x\mapsto x+1$. The group $GA$ is a subgroup of $PSL(2,\R)$.
The group $PSL(2,\R)$ acts on the circle at infinity 
$S^1_\infty$ of the Poincar\'e upper half plane,
where $GA$ is the isotropy subgroup of
$\infty\in S^1_\infty$.
 Cutting 
$S^1_\infty$ at $\infty$, we get a $C^\infty$ action
of $BS(1,2)$ on a compact interval, say $[-1,1]$. 
This is called the {\em affine} action of $BS(1,2)$.

T. Tsuboi \cite{Ts} showed that there is a homeomorphism $h$ of
$[-1,1]$ which is a $C^\infty$ diffeomorphism on $(-1,1)$ such 
that the conjugate by $h$ of any element of $BS(1,2)$ is 
$C^\infty$ tangent to the identity at the end points.
Then the conjugated action
 extends to an $C^\infty$ action on $[-2,2]$ in such a way
that it
is trivial on $[-2,-1]\cup[1,2]$. 
Consider a subset $S^{n-1}\times[-2,2]$ embedded in 
$B$. The group $GA$ acts on $S^{n-1}\times[-2,2]$,
trivially on the first factor.
This way we obtain a subgroup of $\HH(B)$ isomorphic to
$BS(1,2)$, which we shall denote by $G$.

The key fact for the proof of Theorem \ref{T'} is the following
result of \cite{N'}, which improves a semiconjugacy result in \cite{GL}.

\begin{theorem} \label{Navas} {\rm (C. Bonatti, I. Monteverde, A. Navas
 and C. Rivas)} 
Assume $BS(1,2)$ acts faithfully on a compact interval by $C^1$
 diffeomorphisms in such a way that there is no
 interior global fixed poit. Then the action
is topologically conjugate to the  affine action.
In particular, all
the interior orbits are dense.
\end{theorem}

 The compactness assumption on the interval is
indispensible. In fact, there is a $C^\infty$ exotic action of $BS(1,2)$
on $\R$. See \cite{CC}.
Thus in order to apply the above theorem,
we need the following fixed point result in the first place.

\begin{proposition}\label{11}
The fixed point set $\Fix(\Phi(G))$ is nonempty.
\end{proposition}

\bd
We assume for contradiction that $\Fix(\Phi(G))=\emptyset$.
The proof follows the same line as Proposition \ref{p10}. But
since our target manifold is $\R$ and is noncompact, extra care
will be needed. 

Since $\Fix(\Phi(G))=\emptyset$, any orbit of $\Phi(G)$ is unbounded
towards both directions. Since $G$ is finitely generated, $\Phi(G)$
has a compact cross section $I$ in $\R$, that is, a compact interval $I$
which intersects any $\Phi(G)$ orbit. 
In fact, choose any point $x_0\in\R$ and let $x_1$ be the supremum of
$g(x_0)$, where $g$ runs over a finite symmetric generating set.
Then clearly any orbit intersects the interval $I=[x_0,x_1]$.
Since $G\subset\HH(B)$, $I$
is also a cross section for $\Phi(\HH(B))$.
That is, any $\Phi(\HH(B))$ orbit intersects the compact interval $I$.

Now we follow the proof of Proposition 6.1 in \cite{DKNP}, to show that 
there is a unique minimal set $X$ for $\Phi(\HH(B))$.
In fact we shall show a bit more: there is
a nonempty $\Phi(\HH(B))$ invariant closed subset $X$ in $\R$ which
has the property that 
any nonempty $\Phi(\HH(B)))$ invariant closed  subset
contains $X$.

The proof goes as follows.
Let $F$ be the family 
of nonempty $\Phi(\HH(B))$ invariant closed  subsets of $\R$,
and  $F_I$ the family of nonempty closed subsets $Y$ in $I$ such that
 $\Phi(\HH(B))(Y)\cap I=Y$,
where we denote
$$\Phi(\HH(B))(Y)=\bigcup_{g\in\HH(B)}\Phi(g)(Y).$$
Define a map $\phi:F\to F_I$ by $\phi(X)=X\cap I$,
and $\psi:F_I\to F$ by $\psi(Y)=\Phi(\HH(B))(Y)$.
They satisfiy $\psi\circ\phi=\phi\circ\psi={\rm id}$.

Let $\{Y_\alpha\}$ be a totally ordered set in $F_I$.
Then the intersection $\cap_\alpha Y_\alpha$ is nonempty.
Let us show that it belongs to $F_I$, namely,

\begin{equation}\label{e2}
\Phi(\HH(B))(\cap_\alpha Y_\alpha)\cap I=\cap_\alpha Y_\alpha.
\end{equation}
For the inclusion $\subset$, we have
$$
\Phi(\HH(B))(\cap_\alpha Y_\alpha)\cap I\subset 
(\cap_\alpha \Phi(\HH(B))(Y_\alpha))\cap I=
\cap_\alpha (\Phi(\HH(B))(Y_\alpha)\cap I)=\cap_\alpha Y_\alpha.
$$
For the other inclusion, notice that
$$\cap_\alpha Y_\alpha\subset \Phi(\HH(B))(\cap_\alpha Y_\alpha)
\mbox{ and }\cap_\alpha Y_\alpha\subset I.
$$

Therefore by Zorn's lemma, there is a minimal element $Y$ in $F_I$. 
The set $Y$ is not finite. In fact, if it is finite,
the set $X=\psi(Y)$ in $F$ is discrete, and 
there would be a nontrivial homomorphism from $\Phi(\HH(B))$
to $\Z$, contrary to the fact that $\HH(B)$, and hence $\Phi(\HH(B))$,
is simple.

Now the correspondence $\phi$ and $\psi$ preserve the inclusion.
This shows that  
there is no nonempty $\Phi(\HH(B))$
invariant closed proper subset of $X=\psi(Y)$.
In other words, any $\Phi(\HH(B))$ orbit 
contained in
 $X$ is dense in $X$.
Therefore $X$ is either $\R$ itself or a locally
Cantor set. 
In the former case, any nonempty $\Phi(\HH(B))$ invariant
closed subset must be $\R$ itself. 
 
Let us show that in the latter case,
$X$ satisfies the desired property: $X$ is contained in
any nonempty $\Phi(\HH(B))$ invariant closed subset.
For this, we only need to show that
the  $\Phi(\HH(B))$ orbit of any point $x$ in $\R\setminus X$ 
accumulates to a point in $X$. 
In fact, if this is true, then any nonempty $\Phi(\HH(B))$ invariant
closed subset must intersects $X$. But the intersection must be the
whole $X$ by the above remark.

Let $(a,b)$ be the connected component
of $\R\setminus X$ that contains $x$. (If $x\in X$, there is nothing to
prove.)
There is a sequence $g_k\in\HH(B)$ ($k\in\N$) 
such that $\Phi(g_k)(a)$ accumulates
to $a$ and that $\Phi(g_k)(a)$'s are mutually distinct. Then the intervals
$\Phi(g_k)((a,b))$ are mutually disjoint, and consequently
$\Phi(g_k)(x)$ converges to $a$.
This concludes the proof that $X$ is contained in
any nonempty $\Phi(\HH(B)))$ invariant closed subset.

\medskip

Choose $B'\in\BB$ such that $B'\cap B=\emptyset$. 
Any element of $\HH(B')$ commutes with any element of $\HH(B)$.
Define $\FF(B')$ to be the subset of the group $\HH(B')$ 
consisting of those elements
$g$ such that $\Fix(\Phi(g))\neq\emptyset$. 
By a theorem of H\"older, there is a nontrivial element in $\FF(B')$.
For any $g\in\FF(B')$, the set $\Fix(\Phi(g))$ is closed, nonempty
and invariant by $\Phi(\HH(B))$ by the commutativity. 
Therefore we have 
\begin{equation}\label{e3}
X\subset \Fix(\Phi(g)) \mbox{ for any } g\in\FF(B').
\end{equation}
 This shows that
$\FF(B')$ is a subgroup of $\HH(B')$,
normal
and nontrivial. Since $\HH(B')$ is simple, we have
$\FF(B')=\HH(B')$. Finally again by (\ref{e3}), we get
$\Fix(\Phi(\HH(B')))\neq\emptyset$. Since $\HH(B)$ is conjugate
to $\HH(B')$ and $G$ is  a subgroup of $\HH(B)$, we have
$\Fix(\Phi(G))\neq\emptyset$, contrary to the assumption.
The contradiction concludes the proof of Proposition \ref{11}.
\qed

\section{Fixed point set of $\Phi(\HH_B)$}

For $B\in\BB$, define a subgroup $\HH_B$ of $\HH$ by
$$\HH_B=\{g\in\HH\mid g={\rm id}\mbox{ in a neighbourhood of }B\}.$$
Let $\Phi:\HH\to\Diff^1_+(\R)$ be a homomorphism satisfying
Assumption \ref{a}.
The purpose of this section is to show the following.

\begin{proposition}\label{p41}
For any $B\in\BB$, the fixed point set $\Fix(\Phi(\HH_B))$ is nonempty.
\end{proposition}

{\sc Proof}. 
Any element of
$\Phi(\HH(B))$ commutes with any element of $\Phi(\HH_B)$.
Let us denote $F=\Fix(\Phi(G))$, which we have shown to be nonempty
in Proposition \ref{11}. Clearly $F$ is invariant by any element
of $\Phi(\HH_B)$.
We shall show that there is a fixed point of
$\Phi(\HH_B)$ in $F$. If $F$ is bounded to the left or to the
right, then the extremal point will be a fixed point of $\Phi(\HH_B)$.
So we assume that $F$ is unbounded towards both directions. That is,
any connected component $U$ of $\R\setminus F$ is bounded.

Assume that there is $g\in\HH_B$ such that 
$\Phi(g)(U)\cap U=\emptyset$. (Otherwise $\Phi(g)(U)=U$ for any
$g\in\HH_B$,
and the proof 
will be complete.)
There is a subgroup  $G'$ of $\HH_B$ conjugate to $G$.
By some abuse, denote the generators of $G'$ by $a$ and $b$.
They satisfy $aba^{-1}=b^2$.
Notice that finite products 
of conjugates of $b^{\pm1}$ by elements of $\HH_B$
form a normal subgroup of $\HH_B$.
Since $\HH_B$ is simple, any element of $\HH_B$ can be written as such a product.
Writing the above element $g$ this way, one finds a conjugate of $b$ whose
$\Phi$-image displaces $U$. We may
assume
 that
$\Phi(b)U\cap U=\emptyset$, passing from $G'$ to its conjugate by
an element of $\HH_B$ if necessary. (The conjugate is still denoted by $G'$.)

 Let $V$ be the component of
$\R\setminus\Fix(\Phi(G'))$
that contains $U$. Since $G'$ is conjugate to $G$,
$V$ is a bounded open interval
and $F\cap V$ is a  closed nonempty proper subset of $V$
invariant by $\Phi(G')$. It is easy to show that
$\Phi(b)\vert_{V}\neq{\rm id}$ implies that the action
$\Phi(G')\vert_V$ is faithful.
By Theorem \ref{Navas},
any $\Phi(G')$ orbit in $V$ must be dense in $V$.
This contradicts the fact that $F\cap V$ is invariant
by $\Phi(G')$. The proof is now complete.
\qed

\section{Proof of Theorem \ref{T'}}

Again we assume that $\Phi:{\mathcal G}\to{\rm Diff}^1_+(S^1)$ is
a homomorphism satisfying Assumption \ref{a}. Our purpose here is to get
a contradiction. We follow
an argument in \cite{Mi2}.
\begin{lemma}\label{l51}
Assume $B$ and $B'$ are mutually disjoint balls of $\BB$.
Then any $g\in\HH$ can be written as $g=g_1\circ g_2\circ g_3$,
where $g_1$ and $g_3$ belongs to $\HH_B$ and $g_2$ to $\HH_{B'}$.
\end{lemma}

\bd Take any $g\in\HH$. Then there is an element $g_1\in\HH_B$ such
that $g_1^{-1}\circ g(B)$ is disjoint from $B'$. Next, there is an
element
$g_2\in\HH_{B'}$ such that $g_2^{-1}\circ g_1^{-1}\circ g$ is the
identity in
a neighbourhood of $B$. Thus $g_3=g_2^{-1}\circ g_1^{-1}\circ g$ belongs
to
$\HH_B$ and the proof is complete.
\qed

\begin{lemma} \label{l52}
Assume $B$ and $B'$ are mutually disjoint elements of $\BB$. If two
 points $a$ and $b$ ($a<b$) belong to $\Fix(\Phi(\HH_{B}))$, then
$\Fix(\Phi(\HH_{B'}))\cap[a,b]=\emptyset$.
\end{lemma}

\bd 
Assume a point $c$ in $[a,b]$ belongs to $\Fix(\Phi(\HH_{B'}))$.
Choose an arbitrary element $g\in\HH$. There is
a decomposition $g=g_1\circ g_2\circ g_3$ as in
Lemma \ref{l51}. 
Now $\Phi(g_3)(a)=a$. Since $\Phi(g_2)(c)=c$ and $a\leq c$, we have
$\Phi(g_2)\circ \Phi(g_3)(a)\leq c$. Likewise 
$\Phi(g)(a)=\Phi(g_1)\circ \Phi( g_2)\circ \Phi( g_3)(a)\leq b$.
Since $g\in\GG$ is arbitrary, the $\Phi(\GG)$ orbit of
$a$ is bounded from the right. 
Then the supremum of the orbit must be a global fixed poit,
which is against Assumption \ref{a}:
$\Phi(\GG)$ has no global fixed point.
\qed

\medskip
For any point $x\in \R^n$,
define a subgroup $\HH_x$ of $\HH$ by
$$
\HH_x=\{g\in\HH\mid g \mbox{ is the identity in a neighbourhood of }
x\}.$$

\begin{lemma}\label{l4}
For any $x\in\R^n$, the fixed point set $\Fix(\Phi(\HH_x))$ is nonempty.
\end{lemma}

\bd 
Notice that for any $x\in\R^n$, there is an decreasing sequence
$\{B_k\}$ ($k\in\N$) in $\BB$ such that $\{x\}=\bigcap_kB_k$.
Then $\HH_{B_k}$ is an increasing senquence of subgroups
of $\HH$ such that $\bigcup_k\HH_{B_k}=\HH_x$.
Therefore the closed subsets $\Fix(\Phi(\HH_{B_k}))$ is decreasing and
we have $$\Fix(\Phi(\HH_x))=\bigcap_k\Fix(\Phi(\HH_{B_k})).$$

Now it suffices to prove that $\Fix(\Phi(\HH_B))$ is compact for
 $B\in\BB$. 
Recall that we have already shown that  $\Fix(\Phi(\HH_B))$ is nonempty. 
Assume in way of contradiction that
$\sup\Fix(\Phi(\HH_B))=\infty$. (The other case can be dealt with
similarly.)
Choose $B'\in\BB$ such that $B\cap B'=\emptyset$. Notice that
$\Phi(\HH)$ consists of orientation preserving
diffeomorphisms and $\Phi(\HH_{B'})$ is conjugate to
$\Phi(\HH_{B})$ by such a diffeomorphism. Therefore we also have that
$\sup\Fix(\Phi(\HH_{B'}))=\infty$. 
Now one can find points $a,b\in\Fix(\Phi(\HH_B))$ and a point
$c\in\Fix(\Phi(\HH_{B'})$ such that $a<c<b$. This is contrary to
Lemma \ref{l52}. 
\qed

\bigskip
We use the assumption $n\geq2$ only in the sequal.
Let $D_0$ be the unit compact disc centered at $0$ in $\R^{n-1}\subset\R^n$.
Define a family $\DD$ of closed subsets of $\R^n$ by
$$\DD=\{g(D_0)\mid g\in\HH\}.$$  
For any $D\in\DD$, define a subgroup $\HH_D$ of $\HH$ by
$$
\HH_D=\{g\in\HH\mid g \mbox{ is the identity in a neighbourhood of }
D\}.$$

Lemma \ref{l4} implies that $\Fix(\Phi(D))\neq\emptyset$ for any
$D\in\DD$.

\begin{lemma}\label{l5}
For any $D\in\DD$, the set $\Fix(\Phi(D))$ is a singleton.
\end{lemma}

\bd
First of all notice that for any $D,D'\in\DD$ such that $D\cap
D'=\emptyset$,
we have $\Fix(\Phi(\HH_D))\cap\Fix(\Phi(\HH_{D'}))=\emptyset$.
In fact, as is easily shown, $\HH_D$ and $\HH_{D'}$ generate $\HH$.
Thus the point of the above intersection would be a global fixed point
of $\HH$, against Assumption \ref{a}.
This shows that the interior
 $\Int(\Fix(\Phi(\HH_D)))$ is empty.
In fact, there are
uncountably many mutually disjoint elements of $\DD$,
while mutually disjoint open subsets of $\R$
are at most countable.

Assume that $\Fix(\Phi(\HH_D))$ contains more than one points. 
Since $\Int(\Fix(\Phi(\HH_D)))$ is empty, 
 $\Fix(\Phi(\HH_D))$ is not connected.
To any
$D\in\DD$,
assign a bounded component $I_D$
of $\R\setminus\Fix(\Phi(\HH_D))$ in an arbitrary way. This is possible
by the axiom
of choice. Notice that Lemmata \ref{l51}
and \ref{l52} for the family $\BB$ are valid for
 $\DD$ as well.  (No changes of the proofs are needed.)
Consequently  $I_D\cap I_{D'}=\emptyset$
if $D\cap D'=\emptyset$. Again this is contrary to the fact that  
 there are
uncountably many mutually disjoint elements of $\DD$.
\qed

\medskip
Finally let us prove Theorem \ref{T'}. 
Choose any element $D\in\DD$ and distinct two ponits $x_1,x_2\in D$ 
that are contained in $D$. Then since $\Fix(\Phi(\HH_D))$ is a singleton
and $\Fix(\Phi(\HH_{x_i}))$ is nonempty, we have
$\Fix(\Phi(\HH_{x_1}))=\Fix(\Phi(\HH_{x_2}))$. But
 $\HH_{x_1}$ and $\HH_{x_2}$ generate $\HH$, and there would be
 a global fixed point of $\Phi(\HH)$, against Assumption \ref{a}.
The contradiction shows that the homomorphism $\Phi$ must be trivial.

\section{Sporadic results for $\Diff^0(S^1)$ target}

Let $M=L\times S^m$ be a closed $n$-dimensional manifold such that
$1\leq m\leq n$. Then we have the following result.

\begin{theorem}\label{sporadic}
If $n\geq2$ and $r\neq n+1$, there is no nontrivial homomorphism
from $\Diff^r_c(M)_0$ to $\Diff^0(S^1)$.
\end{theorem}

\bd
Assume that $\Phi:\Diff^r_c(M)_0\to\Diff^0(S^1)$ is a nontrivial
homomorphism. The Lie group $PSL(2,\R)$ acts on $S^m$ as
Moebius transformations. 
So it acts on $M=L\times S^m$, trivially on $L$-coordinates.
Denote the inclusion by
$\iota:PSL(2,\R)\to\Diff^r_c(M)_0$.
The simplicity of the group $\Diff^r_c(M)_0$ shows that 
the homomorphism
$$\Phi\circ\iota:PSL(2,\R)\to \Diff^0(S^1)$$
is nontrivial. 

Now Theorem 5.2 in \cite{M''} asserts that the homomorphism
$\Phi\circ\iota$ is the conjugation of
the standard embedding $\iota_0:PSL(2,\R)\to\Diff^0(S^1)$
by a homeomorphism of $S^1$. It is no loss of generality to assume that
$\Phi\circ\iota=\iota_0$, by changing $\Phi$ if necessary.
If the dimension of $L$ is positive, then $\Diff^r_c(L)_0$
also acts on $M$, trivially on $S^m$-coordinates. 
Any element of the group $\Phi(\Diff^r_c(L))$ must commute with
any element of $PSL(2,\R)$. But there is no nontrivial element
in $\Diff^0_+(S^1)$ which commutes with all the element of
$PSL(2,\R)$. A contradiction.

Let us consider the case where $L$ is a singleton. Then there is an element $g$ in
$\Diff^r(S^n)$, $n\geq2$, which commutes with all the elements
of $\iota(SO(2))$ and is not contained in $\iota(SO(2))$.
But $\Phi\circ\iota(SO(2))=SO(2)$, and any
element in $\Diff^0(S^1)$ which commutes with all the elements of $SO(2)$
must be an element of $SO(2)$, contradicting the injectivity
of $\Phi$.

\begin{remark}
There are a wider class of manifolds for which the above argument holds.
For example, if $M$ is the unit tangent bundle of a closed hyperbolic
 surface and if $r\neq4$, then any homomorphism from $\Diff^r_c(M)_0$
to $\Diff^0(S^1)$ is  trivial. 
\end{remark}

\end{document}